\documentclass[12pt]{amsart}
\usepackage{amssymb}
\usepackage{amsmath}
\newtheorem{prop}{Proposition}[section]
\newtheorem{thm}{Theorem}

\newtheorem{definition}[prop]{Definition}  
  
\newtheorem{ex}[prop]{Example}

\def\R{\Bbb R} \def\N{\Bbb N} \def\Z{\Bbb Z} 
 
  \def\Q{\Bbb Q}  

\title{Minimal Homeomorphisms on low-dimension tori }
 
\author{N. M. dos Santos$\dag$} 
\address{Universidade Federal Fluminense, 24020-005 Niteroi, R.J. Brazil}
\email{nathan@mat.uff.br}
\thanks{$\dag$ Corresponding address: Rua Lopes Quintas, 225 ap. 401-A, Jardim Bot\^anico, Rio de Janeiro, Cep 22460-010. Brazil}

\author{R. Urz\'{u}a-Luz$\ddag$}
\address{Universidad Cat{\'o}lica de Norte, Casilla 1280, Antofagasta, Chile}
\email{rurzua@ucn.cl}
\thanks{$\ddag$ Partially supported by Chilean FONDECYT Grant N.1060977}

\date{}

\begin{document}
\maketitle

\begin{abstract}
In this article we study minimal homeomorphisms(all orbits are dense) of the tori $T^{n},$ $n<5.$ The linear part of a homeomorphism %%@
$\varphi $ of $T^{n}$ is the linear mapping $L$ induced by $\varphi $ on the first homology group of $T^{n}$. It follows from the %%@
Lefschetz fixed point theorem that $1$ is an eigenvalue of $L$ if $\varphi $ minimal. We show that if $\varphi $ is minimal and $n<5$ %%@
then $L$ is quasi-unipontent, i.e., all the eigenvalues of $L$ are roots of unity and conversely if $L\in GL( n,\Z)$ is quasi-unipotent %%@
and $1$ is an eigenvalue of $L$ then there exists a $ C^{\infty }$ minimal skew-product diffeomorphism $\varphi $ of $T^{n}$ whose linear %%@
part is precisely $L.$ We do not know if these results are true for $n>4$. We give a sufficient condition for a smooth skew-product %%@
diffeomorphism of a torus of arbitrary dimension to be smoothly conjugate to an affine transformation. 
\end{abstract}

\maketitle 

\section{Minimal homeomorphisms on low-dimension tori}

We first prove
\begin{prop}
Let $\varphi$ be a minimal homeomorphism of a torus $T^{n}$ and $L$ be the induced mapping on $H_{1}( T^{n},\Z)$. Then the minimal %%@
polynomial $p(x)$ of $L$ can not be decomposed over $\Q[x]$, as $p(x)=q(x)r(x)$ where all the roots of $q(x)$ are roots of unity and %%@
$r(x)$ is not constant with no roots in the unit circle.
\end{prop}

\begin{proof}
Assume that $p(x)$ has such a decomposition. Then by the Primary Decomposition Theorem we have an invariant direct sum decomposition over %%@
$\Q$ 
\begin{equation}
\R^{n}=E\oplus V
\label{01}
\end{equation}
where the restriction $B$ of $L$ to $V$ is hyperbolic. Now $\Gamma =V\cap \Z^{n}$ is a discrete cocompact subgroup of $V$ and 
$M=V/\Gamma$ is homeomorphic to a torus $T^{k}$, $k<n$.

Let $b$ be the hyperbolic diffeomorphism of $M$ induced by $B$ and $\varphi$ be given on the covering by $L+F$, where %%@
$F:\R^{n}\rightarrow \R^{n}$ is continuous and $F(x+\ell ) =F(x) $ for all $x\in \R^{n}$ and $\ell \in \Z^{n}$. We claim that $b$ is a
factor of $\varphi$. For, consider the continuous surjective mapping $ h:T^{n}\rightarrow M$ given on the covering $\R^{n}$ by
\begin{equation}
h(x) =P(x)+H(x)
\end{equation}
where $P:\R^{n}\rightarrow V$ is the projection associated to the decomposition $(\ref{01}) $ and $H:T^{n}\rightarrow V$ is a continuous %%@
solution of
the cohomological equation 
\begin{equation}
BH(x) -H(\varphi(x)) =P(F(x))\label{02}
\end{equation}
Now since $B$ is hyperbolic then a continuous solution of $(\ref{02})$ exists see ~\cite{K-H} [\textit{Theorem 2.9.2}] and since $P\circ %%@
L=B\circ P$ then $h\circ \varphi =b\circ h$. Observing that $h\circ \varphi ^{\ell }=b^{\ell}\circ h$ for all $\ell \in \Z$ and since $h$ %%@
is surjective we see that $\varphi$ can not be minimal because $b$ has periodic points.
$\hfill \square$
\end{proof}

\begin{thm}
Any minimal homeomorphism $\varphi $ of a torus $T^n$,  $n<5$ is quasi-unipotent on the homology and $1$ is an eigenvalue of its linear %%@
part.
\end{thm}

\begin{proof}
Minimality of $\varphi $ and the Lefschetz fixed point Theorem shows that $1$ is a root of the minimal polynomial $p(x)$ of the linear %%@
part $L$
of $\varphi$.  Thus $p( x)=(x-1)s(x)$ where $\deg s( x) <4$, since $\deg p<5$.

If $\deg s(x) =3$ then $s(x)$ factors over $\Z\left[ x\right]$ as $(x\pm 1) q( x)$ and by Proposition $1$ all the roots of $q$ are roots %%@
of unity. If $\deg s( x) <3$ again by Proposition $1$ all the roots of $s(x)$ are roots of unity.
$\hfill \square$
\end{proof}

We do not know if the above Theorem is true if $n>4$. There are irreducible polynomials in $\Q\left[ x\right]$ with roots of absolute %%@
value $1$ and roots of absolute value different of $1$.

\begin{ex}
Eisenstein's criterion shows that the polynomial $p(x)=x^{4}+4x^{3}-6x^{2}+4x+1$ is irreducible over $\Q\left[ x\right] $
and as 
\begin{equation*}
p(x) =( x^{2}+2(1-\sqrt{3})x+1)(x^{2}+2( 1+\sqrt{3}) x+1)
\end{equation*}
we can see that 
\begin{equation*}
(\sqrt{3}-1) \pm i\lambda ,\text{ }\lambda =\sqrt{(1-( \sqrt{3}-1) ^{2}) }
\end{equation*}
are roots of absolute value $1$ and they are not roots of unity and the other two roots of $p(x)$ have absolute value different from $1$.
\end{ex}

\section{Minimal skew-product transformations of the torus.}

In this section we show that every quasi-unipotent matrix $L\in GL( n,\Z)$, $n<5$ with $1$ as eigenvalue is the linear part of a smooth %%@
minimal skew-product transformation of the torus $T^{n}$. Actually the skew - products are of the
particular type given in $(\ref{2})$. Notice that $T^{n-p}$ acts freely on $T^{p}\times T^{n-p}$ by translation on the second
factor. Thus if a homeomorphism $\psi $ of $T^{n}$ commutes with this action it induces a homeomorphism $\psi _{0}$ of the orbit space 
$T^{n}/T^{n-p}$ which is homeomorphic to the torus $T^p$ and we say that $( T^n,\psi )$ is a free $T^{n-p}$-extension of $(T^p,\psi _{0}) %%@
$  ~\cite{P}
\begin{thm}
Let $L\in GL( n,\Z) ,$ $n<5$ be quasi-unipotent having $1$ as an eigenvalue. Then there exists a minimal smooth skew-product %%@
diffeomorphism $\varphi $ of the torus $T^{n}$ whose linear part is $L$.
\end{thm}

\begin{proof}
We may assume that by $\left[ \text{Newman}\right]$ ~\cite{N}
\begin{equation}
L=\left( 
\begin{array}{cc}
A & \mathbf{0} \\ 
C & B
\end{array}
\right)  \label{1}
\end{equation}
where $A\in GL( p,\Z) ,$ $A=I+N_{1}$ and $B \in GL( n-p,\Z) $ such that $B^{m}=I+N_{2}$, 
$m\in \Z^+$,  where $N_{1}$ and $N_{2}$ are nilpotent.

Let $\varphi $ be diffeomorphism of $T^{n}$ given on the covering $\R^n$ by
\begin{equation}
\varphi ( X,Y) =( AX+\alpha ,CX+BY+F( X) )
\label{2}
\end{equation}
where $a(X) =AX+\alpha $ gives an affine minimal transformation of $T^{p}$ ~\cite{H} and $F:\R^{p}\rightarrow \R^{n-p}$ is a smooth %%@
$\Z^{p}-$periodic function, i.e., $F(X+\ell) =F(X)$ for all $\ell \in \Z^{p}$

The iterates of $\varphi $ are given by
\begin{equation}
\varphi ^{m}( X,Y) =( a^{m}( X)
,C_{m}X+B^{m}Y+\alpha ( m) +F_{m}( X) )
\label{3}
\end{equation}
where
\begin{equation*}
C_{l}X=\sum_{j=1}^{l}B^{k-j}CA^{j-1}
\end{equation*}
\begin{equation*}
\alpha ( m) =( \sum_{j=1}^{m}C_{m-j}) \alpha
\end{equation*}
 \begin{equation}
 F_{m}( X) =\sum_{j=1}^{m}B^{m-j}F( a^{j-1}( X))
\label{4}
\end{equation}

Recall that an action is \textit{simple} if for each character $\gamma \in \hat{T}^{n-p}$ there exist a continuous function %%@
$f_\gamma:T^{n}\rightarrow T^{1}$ and $g:T^p\rightarrow T^1$ such that
\begin{equation}
f_{\gamma }(z,w)=g(z)\gamma (w)
\label{5}
\end{equation}
for every $z\in T^p$ and $w\in T^{n-p}$.

Let $\nu $ be the endomorphism of $T^{n-p}$ given on the covering $\R^{n-p}$ by the nilpotent matrix $N_{2}$ and consider the subgroup %%@
$H=\ker \nu$ of $T^{n-p}$. Notice that $\psi =\varphi ^{m}$ is invariant under the restriction of the action of $T^{n-p}$ to the subgroup %%@
$H$, i.e., 
\begin{equation}
\psi ( h\cdot (z,w) ) =h\cdot \psi ( z,w)
\label{6}
\end{equation}
$h\in H$

From now on we will assume that the diffeomorphism $\psi _{0}$ of $T^{n}/H$ induced by $\psi $ is either minimal or uniquely ergodic.

Let $\pi :T^{n}\rightarrow T^{n}/H$ be the projection . By $\left[ \text{Parry, Theorem 1}\right] $ ~\cite{P} $\psi $ is minimal(uniquely %%@
ergodic) if only if
the equation 
\begin{equation}
\frac{f(\psi _{0}( \pi ( z,w) ) ) }{f( z,w) }=\frac{f_{\gamma }( \psi (z,w) ) }{f_{\gamma }( z,w) }  \label{7}
\end{equation}
has no continuous(measurable) solution for each character $\gamma \in \hat{H}$, $\gamma \neq 1$. Observe that condition in $(\ref{7})$ %%@
does not depend on the choice of the particular function $f_{\gamma }$.
Notice that the functions $f:T^{n}/H\rightarrow T^{1}$ are given by the $H$-invariant functions $f:T^{n}\rightarrow T^{1}$, i.e., 
$f(z,hw)=f(z,w)$ for all $h\in H.$

If $L$ is unipotent then there exists a minimal affine diffeomorphism of $T^{n}$ ~\cite{H}. Thus it suffices consider $L$ quasi-unipotent %%@
but not unipotent.

Suppose $n=2.$ We may assume ~\cite{N}
\begin{equation}
L=\left( \begin{array}{cc}
1 & 0 \\ 
s & -1
\end{array}
\right)  \label{8}
\end{equation}
Let $\varphi ( x,y) =( x+\alpha ,sx-y+F( x))$. Thus 
\begin{equation*}
\psi ( x,y) =\varphi ^{2}( x,y) =( x+2\alpha,s\alpha +y-F(x) +F( x+\alpha))
\end{equation*}
It is easy to see that $( T^{2},\psi )$ is simple free $T^{1}$-extension of the translation of $T^{1}$ given on the covering $\R^{1}$ by 
\begin{equation*}
\psi _{0}(x) =x+2\alpha.
\end{equation*}
We choose a Liouville number $\alpha$ and a sequence $\left\{ k_{j}\right\} _{j\in N}$ so that 
\begin{equation*}
|e^{2\pi ik_{j}\alpha}+1|<\frac{4\pi}{( k_{j}) ^{j}}
\end{equation*}
 see [Appendix] and $F:T^{1}\rightarrow \R$ given by the Fourier transform 
\begin{equation}
\hat{F}( k) =\left\{ 
\begin{array}{cc}
0 & \text{if }k\neq \pm k_{j} \\ 
1+e^{\pm 2\pi ik_{j}\alpha } & \text{ if }k=\pm k_{j}
\end{array}
\right.  \label{9}
\end{equation}
By $(\ref{7}) $ it suffices to show that the cohomological equation 
\begin{equation}
\frac{f( \psi _{0}(x)) }{f( x) }=\frac{f_{\gamma }( \psi ( x,y) ) }{f_\gamma(x,y)}  \label{10}
\end{equation}
where $f$ is given on the covering $\R^{1}$ by $f( x) =\ell x+G( x)$ and $f_\gamma$ is given on the covering $\R^{2}$ by 
$f_\gamma( x,y) =\ell _{\gamma }y$, $\ell$ and $\ell _\gamma$ in $\Z$, has no continuous solution $f:T^{1}\rightarrow \R^{1}$. 
This is equivalent to show that the equation 
\begin{equation}
\ell \alpha +G(x+2\alpha )-G( x) =\ell _\gamma s\alpha +\ell _\gamma\left[ F(x+\alpha)-F(x)\right]  \label{11}
\end{equation}
has no continuous solution $G$. If $\ell\neq\ell _\gamma$ then one see that $(\ref{11})$ has no continuous solution $G$. 
If $\ell =s\ell _{\gamma }$ the equation $(\ref{11})$ becomes
\begin{equation*}
G(x+2\alpha ) -G( x) =\ell _{\gamma }\left[ F(x+\alpha ) -F( x) \right]
\end{equation*}
and gives the Fourier coefficients equations 
\begin{equation}
\hat{G}(k) =\left\{ 
\begin{array}{cc}
0 & \text{if }k\neq \pm k_{j} \\ 
\ell _\gamma & \text{ if }k=\pm k_{j}
\end{array}
\right.
\end{equation}

which by choice of $F$ does not give a $L^{1}-$solution $G$.

Suppose now $n=3$. There are two possibilities for the characteristic polynomial $p(x)$ of $L$, $p( x)=(x-1)q( x)$ or 
$p(x) =(x-1)^2(x+1)$ where $q(1)\neq 0$.

If $p( x) =( x-1) q( x)$ then we assume that ~\cite{N} 
\begin{equation*}
L=\left( 
\begin{array}{cc}
1 & 0 \\ 
C & B
\end{array}
\right)
\end{equation*}
where $B\in GL( 2,\Z)$ is quasi-unipotent and $1$ is not an eigenvalue of $B$. Thus either $B$ is periodic with period $m=3,4,6$ or 
$B=\left( \begin{array}{cc}
-1 & 0 \\ 
s & -1
\end{array}\right)$, $s\in Z$ and by $( \ref{2})$ the diffeomorphism $\varphi$ is given by 
\begin{equation}
\varphi (x,Y) =(x+\alpha ,Cx+BY+F(x))
\label{-1}
\end{equation}
where $Y=(y,z)$ and $F:\R\rightarrow \R^{2}$ is a $\Z$-periodic function.

If $B$ periodic then by $(\ref{3})$
\begin{equation*}
\psi (x,Y) =\varphi ^{m}( x,Y) =( x+m\alpha, Y+\alpha ( m) +F_{m}(x))
\end{equation*}
since $C_{m}=0$.

Hence $( T^{3},\psi)$ is a simple free $T^2$-extension of the translation of $T^{1}$ given on the covering $\R^{1}$ by
\begin{equation*}
\psi _{0}( x) =x+\alpha
\end{equation*}
Choose by [Appendix] a Liouville number $\alpha$ and a sequence of integers $\left\{k_{j}\right\} _{j\in N}$ such that 
\begin{equation}
|e^{2\pi ik_{j}\alpha}-e^{\frac{2\pi i}{m}}|<\frac{4\pi}{( k_{j}) ^{j}}.
\label{13}
\end{equation}
By ~\cite{P} $\psi $ is a minimal diffeomorphism if  the equation $(\ref{7})$ has no continuous solution 
$f(x) =e^{2\pi i\left[ \ell x+G(x)\right]}$ where 
$f_{\gamma }(x,y,z) =e^{2\pi i\left\langle \mathbf{\ell}_\gamma ,( y,z) \right\rangle }$ or equivalently the equation 
\begin{equation}
m\ell \alpha +G( x+\alpha ) -G( x) =\left\langle \mathbf{\ell }_\gamma,\alpha (m)+F_{m}(x)\right\rangle
\label{14}
\end{equation}
has no continuous solution $G$. If $m\ell \alpha \neq \left\langle \mathbf{\ell }_{\gamma },\alpha ( m) \right\rangle$ 
then one see that $(\ref{14})$ has no continuous solutions $G$.  If 
$m\ell \alpha=\left\langle \mathbf{\ell }_{\gamma },\alpha (m) \right\rangle$ the equation $(\ref{14})$ becomes
\begin{equation}
G( x+m\alpha )-G( x) =\left\langle \mathbf{\ell }_{\gamma },F_{m}(x)\right\rangle
\label{c15}
\end{equation}
or in Fourier Coefficients
\begin{equation}
\hat{G}( k) =\left\langle \mathbf{\ell }_{\gamma },( e^{2\pi
ik\alpha }I-B) ^{-1}\hat{F}(k)\right\rangle  \label{15}
\end{equation}
Consider $F:T^{1}\rightarrow \R^{2}$ the smooth function given by the Fourier transform
\begin{equation}
\hat{F}(k)=\left\{
\begin{array}{cc}
0 & k\neq \pm k_{j} \\ 
( e^{2\pi ik_{j}\alpha }-e^{\frac{2\pi i}{m}}) V & k=k_{j} \\ 
( e^{-2\pi ik_{j}\alpha }-e^{-\frac{2\pi i}{m}}) \overline{V} & 
k=-k_{j}
\end{array}
\right.  \label{15'}
\end{equation}
where $V$ is the eigenvector of $B$ associate to the eigenvalue $e^{\frac{2\pi i}{m}}$.  Then 
\begin{eqnarray*}
\hat{G}( k_{j}) &=&\left\langle \mathbf{\ell }_{\gamma },(e^{2\pi ik_{j}\alpha }I-B)^{-1}\hat{F}(k_{j})\right\rangle=\left\langle %%@
\mathbf{\ell }_{\gamma },V\right\rangle \\
\hat{G}( -k_{j}) &=&\left\langle \mathbf{\ell }_{\gamma },(e^{-2\pi ik_{j}\alpha }I-B)^{-1}\hat{F}(-k_{j})\right\rangle=\left\langle %%@
\mathbf{\ell }_{\gamma },\overline{V}\right\rangle
\end{eqnarray*}
As the period of $B$ is $m>2$ then $\left\langle \mathbf{\ell }_\gamma, V\right\rangle \neq 0$ for all 
$\mathbf{\ell}_\gamma\in \Z^{2}-\{ 0\}$. Thus by choice of  $F$ the equation $(\ref{c15})$ has no $L^{1}$-solution.

If $B=\left( 
\begin{array}{cc}
-1 & 0 \\ 
s & -1
\end{array}
\right)$, $s\in \Z$ then by $( \ref{-1}) $
\begin{equation}
\varphi (x,Y) =( x+\alpha ,Cx+BY+F( x) )
\label{14''}
\end{equation}
Thus 
\begin{equation*}
\psi( x,y,z) =\varphi ^{2}( x,y,z) =( x+2\alpha, C_{2}x+\alpha( 2) +B^{2}Y+F_2(x))
\end{equation*}
If $s=0$, $( T^{3},\psi )$ is a simple free $T^{2}$-extension of the translation $(T_{1},\psi _{0})$, given in covering $\R^{1} $ by 
$\psi _{0}(x) =x+2\alpha$.

By ~\cite{P} $\psi$ is a minimal diffeomorphism if the equation $(\ref{7})$ has no continuous solution $f(x) =e^{2\pi i\left[ \ell x+G( %%@
x) \right]}$, where $f_{\gamma }(x,y,z) =e^{2\pi i\left\langle \mathbf{\ell }_{\gamma },( y,z) \right\rangle }$ or equivalently the %%@
equation
\begin{equation}
m\ell \alpha +G( x+2\alpha ) -G( x) =\left\langle \mathbf{\ell }_\gamma ,\alpha ( 2) +F_{2}(x)\right\rangle
\label{14'}
\end{equation}
has no continuous solution $G$

If $2\ell \alpha \neq \left\langle \mathbf{\ell }_{\gamma },\alpha (2) \right\rangle $ then one see that 
$( \ref{14'}) $ has no continuous solutions $G$.
 If $2\ell \alpha =\left\langle \mathbf{\ell }_{\gamma },\alpha ( 2) \right\rangle $ the equation becomes 
\begin{equation*}
G( x+2\alpha )-G( x) =\left\langle \mathbf{\ell }_{\gamma },F_{2}(x)\right\rangle
\end{equation*}
or in Fourier Coefficients
\begin{equation*}
\hat{G}( k) =\left\langle \mathbf{\ell }_{\gamma },( e^{2\pi ik\alpha }+1) ^{-1}\hat{F}(k)\right\rangle
\end{equation*}
Consider $F:T^{1}\rightarrow \R^{2}$ the smooth function given by the Fourier transform 
\begin{equation*}
\hat{F}(k)=\left\{\begin{array}{cc}\mathbf{0} & k\neq \pm k_{j} \\ 
( e^{\pm 2\pi ik_{j}\alpha }+1) V & k=\pm k_{j}
\end{array}
\right.
\end{equation*}
If the vector $V=( a,b) \in \R^{2}$ and $a$ and $b$ are linearly independent over rational numbers, then by the choice of $F$ 
the equation $(\ref{14'})$ has no $L^{1}$-solution.

If $s\neq 0$, consider the diffeomorphim
\begin{equation*}
\varphi ( x,y,z) =( x+\alpha ,px-y+F_{1}( x)
,qx+sy-z)
\end{equation*}
such that $\varphi _0(x,y) =( x+\alpha ,px-y+F_{1}(x))$ is minimal. 
Hence $(T^{3},\psi =\varphi ^{2})$ is a simple free $T^{1}$-extension of the diffeomorphism $(T^{2}, \psi _0)$,  given on covering %%@
$\R^{2}$ by $\varphi _{0}^{2}(x,y) =\psi _{0}( x,y) =( x+2\alpha ,y+p\alpha-F_{1}( x) +F_{1}( x+\alpha))$. By ~\cite{P} $\psi$ is a %%@
minimal diffeomorphism if the equation $( \ref{7})$ has no continuous solution $f( x,y) =e^{2\pi i\left[\left\langle \ell ,( x,y) %%@
\right\rangle +G( x,y)\right]}$, where $f_{\gamma }( x,y,z) =e^{2\pi i\mathbf{\ell}_{\gamma },z}$ or equivalently the equation 
\begin{multline*}
\left\langle \ell ,( 2\alpha ,p\alpha +F_{1}( x+\alpha )
-F_{1}( x) ) \right\rangle +G( \psi _{0}(
x,y) ) -G( x,y)\\
=\ell _{\gamma }\left[
spx-2sy+q\alpha +sF_{1}( x) \right]
\end{multline*}
has no continuous solution $G$.  This is so since the right hand side of the above equation is not a periodic function for $\ell _\gamma$ %%@
since $s\neq 0$.

If $p(x) =( x-1) ^{2}( x+1)$ then we assume that ~\cite{N} 
\begin{equation*}
L=\left(
\begin{array}{cc}
A & 0\\
C & -1
\end{array}
\right)
=\left( 
\begin{array}{ccc}
1 & 0 & 0 \\ 
p & 1 & 0 \\ 
q & s & -1
\end{array}
\right)
\end{equation*}
by $(\ref{2}) $ the diffeomorphism $\varphi $ is given by
\begin{equation*}
\varphi( X,z) =(AX+\delta, CX-z+F(x))=(a(X),CX-z+F(x))
\end{equation*}
where $X=(x,y)$ and $\delta=(\alpha, \beta)$ and $F$ a periodic smooth function, then 
\begin{equation*}
\psi( X,z) =\varphi ^{2}( X,z) =( a^2(X),C_2X+z+\alpha(2) +F_2(x))
\end{equation*}

Hence $(T^3,\psi)$ is a simple free $T^{1}$-extension of the minimal affine transformation $( T^2,\psi _0)$, given on covering $\R^2$ by 
$\psi _0(X) =a^2(X)$. Again, by ~\cite{P} $\psi$ is a minimal diffeomorphism if the equation $(\ref{7})$ has no continuous solution %%@
$f(X)=e^{2\pi i\left[ \left\langle \ell,X \right\rangle +G(X) \right]}$, where $f_\gamma(X,z) =e^{2\pi i\mathbf{\ell }_{\gamma }z}$ or %%@
equivalently the equation
\begin{equation}
\left\langle \mathbf{\ell },a^2(X)-X\right\rangle +G(\psi _0(X))-G(X)=\ell _\gamma \left[C_2X+ \alpha(2)+F_2(x) \right]  
\label{eq1}
\end{equation}
has no continuous solution $G$. If $\left\langle \mathbf{\ell },a^2(X)-X\right\rangle \neq \ell _\gamma \left[C_2X+ \alpha(2)\right]$ %%@
then one see that $(\ref{eq1})$ has no continuous solutions $G$.  \newline If $\left\langle \mathbf{\ell },a^2(X)-X\right\rangle=\ell %%@
_\gamma \left[C_2X+ 
\alpha(2)\right]$ the equation $(\ref{eq1})$ becomes
\begin{equation}
G(\psi _{0}(X)) -G(X) =\ell _\gamma F_2(x) 
\label{eq2}
\end{equation}

If $F:T^{1}\rightarrow \R$ is given in Fourier coefficients as in $(\ref{9})$ the equation $(\ref{eq2})$ has no continuous solution $G$.

 Finally suppose that $n=4$. There are three possibilities for the characteristic polynomial $p(x)$ of $L$, 
$p(x)=( x-1)q_{1}(x)$, $p(x)=(x-1)^2 q_2(x)$ or $p(x) =(x-1)^3( x+1)$, where $q_1$ and $q_2$ are irreducible over $\Q(x)$. 

If $p(x) =(x-1) q_{1}(x)$ then we assume that ~\cite{N}
\begin{equation*}
L=\left( 
\begin{array}{cc}
1 & 0 \\ 
C & B
\end{array}
\right)
\end{equation*}
where $B\in GL( 3,\Z) $ is quasi-unipotent and $1$ is not an eigenvalue of $B$. Thus
\begin{equation*}
B=\left( 
\begin{array}{cc}
-1 & 0 \\ 
C_{0} & B_{0}
\end{array}
\right)
\end{equation*}
where $B_{0}\in GL( 2,\Z) $ is quasi-unipotent by $(\ref{2})$ the diffeomorphism $\varphi$ is given by 
\begin{equation}
\varphi(x,Y) =(x+\alpha ,Cx+BY+F( x)))
\label{16'}
\end{equation}
where $Y=( y,z,w)$ and $F:\R\rightarrow \R^{3}$, $F(x)=(F_1(x) ,F_2(x) ,F_3(x))$ is $\Z$-periodic smooth function.

If $-1$ is not an eigenvalue of $B_0$ then $B_0$ is periodic with period $m=3,4$ and $6$, thus $B^{2m}=I$, then by $(\ref{3})$ we have
\begin{equation}
\psi (x,Y) =\varphi ^{2m}( x,y,Y) =(x+2m\alpha, C_{2m}x+\alpha ( 2m) +Y+F_{2m}(x))  \label{16}
\end{equation}

Choose by [Appendix] a Liouville number $\alpha $ and two sequences of integers $\left\{ k_{j}\right\} _{j\in N}$ and $\left\{ %%@
k_{j}^{\prime }\right\}_{j\in N}$ such that

\begin{equation*}
|e^{2\pi ik_{j}\alpha }-e^{\frac{2\pi i}{m}}|<\frac{4\pi}{( k_{j})^{j}}\text{ and }|e^{2\pi ik_{j}^{\prime }\alpha %%@
}+1|<\frac{4\pi}{(k_{j}^{\prime }) ^{j}}
\end{equation*}
Again, by ~\cite{N}.  It is easy to see that $(T^4,\psi)$ is simple free $T^3$- extension of the translation of $T^1$ given on the %%@
covering $\R^1$ by 
\begin{equation*}
\psi _{0}(x) =x+2m\alpha.
\end{equation*}
Now $\psi$ is a minimal diffeomorphism if the equation $(\ref{7})$ has no continuous solution $f(x) =e^{2\pi i\left[\ell x+G(x)\right]}$ 
where $f_\gamma(x,Y)=e^{2\pi i\left\langle \mathbf{\ell}_\gamma,(y,Y)\right\rangle }$ or equivalently the equation below has no %%@
continuous solution $G$
\begin{equation}
2m\ell \alpha +G(x+2m\alpha ) -G( x) =\left\langle 
\mathbf{\ell }_{\gamma },C_{2m}x+\alpha (2m)
+F_{2m}(x)\right\rangle  \label{17}
\end{equation}
If $2m\ell \alpha \neq \left\langle \mathbf{\ell} _\gamma,\alpha(2m)+ C_{2m}x \right\rangle$ then one see that $(\ref{17})$ has no %%@
continuous solutions $G$. If $2m\ell \alpha =\left\langle \mathbf{\ell} _\gamma,\alpha(2m)+ C_{2m}x \right\rangle$ becomes
\begin{equation}
G( x+2m\alpha ) -G( x) =\left\langle \mathbf{\ell }
_{\gamma },F_{2m}(x)\right\rangle  \label{17'}
\end{equation}
or in Fourier coefficients
\begin{equation*}
\hat{G}(k)=\left\langle \mathbf{\ell}_{\gamma },( e^{2\pi ik\alpha }I-B) ^{-1}\hat{F}(k)\right\rangle
\end{equation*}
Consider $F_1:\R^1\rightarrow \R^1$ smooth function given by the Fourier transform
\begin{equation*}
\hat{F}_{1}(k)=\left\{\begin{array}{cc}
0 & k\neq \pm k_{j}^{\prime } \\ 
( e^{\pm 2\pi ik_{j}^{\prime }\alpha }+1) & k=\pm k_{j}^{\prime }
\end{array}
\right.
\end{equation*}
as in $(\ref{9})$ and the smooth function $F:\R^{1}\rightarrow \R^{2}$,  
$F( x) =( F_{2}( x) ,F_{3}( x) )$ given by the Fourier transform 
\begin{equation*}
\hat{F}(k)=\left\{ 
\begin{array}{cc}
\mathbf{0} & k\neq \pm k_{j} \\ 
( e^{2\pi ik_{j}\alpha }-e^{\frac{2\pi i}{m}}) V & k=k_{j} \\ 
( e^{-2\pi ik_{j}\alpha }-e^{-\frac{2\pi i}{m}}) \overline{V} & k=-k_{j}
\end{array}
\right.
\end{equation*}
as in $( \ref{15'})$ the equation $(\ref{17'})$ has no continuous  solution $G$.

Now if $B_{0}=\left( 
\begin{array}{cc}
-1 & 0 \\ 
s & -1
\end{array}
\right) $
 then
\begin{equation*}
L=\left( 
\begin{array}{cc}
A & \mathbf{0} \\ 
C_{0} & B_{0}
\end{array}
\right) 
\end{equation*}
where 
$A=\left( 
\begin{array}{cc}
1 & 0 \\ 
p & -1
\end{array}
\right)$ and the diffeomorphism given in $(\ref{16'})$ we can written as
\begin{equation*}
\varphi ( X,Y) =(AX+F,C_0X+B_0Y+H(x)) 
\end{equation*}
where $X=(x,y)$, $Y=(z,w)$, $F(x)=(\alpha ,F_{1}(x))$ and $H(x) =(F_2(x), F_3(x))$ are smooth $\Z$-periodic functions. These functions %%@
are determined by the Fourier transforms

\begin{equation}
\hat{F_1}( k) =\left\{ 
\begin{array}{cc}
0 & \text{if }k\neq \pm k_{j} \\ 
1+e^{\pm 2\pi ik_{j}\alpha } & \text{ if }k=\pm k_{j}
\end{array}\label{1'}
\right.
\end{equation}
and 
\begin{equation}
\hat{H}( k) =\left\{ 
\begin{array}{cc}
0 & \text{if }k\neq \pm k_{j} \\ 
(1+e^{\pm 2\pi ik_{j}\alpha })V & \text{ if }k=\pm k_{j}
\end{array}
\right.\label{2'}
\end{equation}
where $V=(a,b)$. Let us consider the diffeomorphism 
\begin{multline*}
\psi ( X,Y) =\varphi ^{2}( X,Y) =(X+AF(x) +F(x+\alpha),\\
C_{0}(2) X+B_0^2Y+C_0 F(x) +B_0H(x) +H(x+\alpha))
\end{multline*}
and suppose that $\varphi _{0}( x,y) =( x+\alpha, px-y+F_{1}( x) ) $ is minimal. Hence \newline $( T^{4},\psi=\varphi ^{2}) $ is a simple %%@
free $T^{2}$-extension of the minimal diffeomorphism $(T^2,\psi _0)$, given on covering $\R^2$ by 
$\varphi _0^2( x,y) =\psi _0(x,y) =(x+2\alpha,y+p\alpha-F_1( x) +F_1( x+\alpha))$. By ~\cite{P} $\psi$ is a minimal diffeomorphism if the %%@
equation $( \ref{7})$ has no continuous solution $f(X) =e^{2\pi i\left[ \left\langle \mathbf{\ell},X\right\rangle
+G( X) \right]}$,  where $f_\gamma(X,Y)=e^{2\pi i\left\langle \mathbf{\ell}_\gamma, Y\right\rangle}$. This is equivalently to the %%@
equation 
\begin{multline}
\left\langle \mathbf{\ell },AF( x+\alpha)+F( x) \right\rangle +G( \psi _{0}( X))-G( X)\\
=\left\langle \mathbf{\ell }_\gamma, C_{0}(2)X+[B_{0}^{2}-I]Y+C_{0}F( x)+B_{0}H( x)+H(x+\alpha)\right\rangle 
\label{16''}
\end{multline}
If $\left\langle \mathbf{\ell }_{\gamma },C_{0}( 2)X+[B_{0}^{2}-I]Y\right\rangle \neq 0$ then one see that
$(\ref{16''})$ has no continuous solution $G$. 
If $\left\langle \mathbf{\ell }_\gamma,C_{0}(2)X+[B_{0}^{2}-I]Y\right\rangle=0$ the equation $( \ref{16''})$ becomes  

\begin{multline}
\left\langle \mathbf{\ell },AF( x) +F( x+\alpha) \right\rangle +G( \psi _{0}( X) )-G( X)\\
=\left\langle \mathbf{\ell}_\gamma, C_{0}F( x)+B_{0}H( x)+H(x+\alpha)\right\rangle 
\label{3'}
\end{multline}
Integranting $(\ref{3'})$ along the fibres of the bundle $(x,y)\rightarrow y$ we get 
\begin{multline}
\langle \mathbf{\ell },AF(x) +F(x+\alpha)\rangle +g( x+\alpha)-g(x)\\
=\left\langle \mathbf{\ell }_\gamma, C_{0}F( x)+B_{0}H( x)+H(x+\alpha)\right\rangle 
\label{16'''}
\end{multline}
where $g(x)=\int _{T^1}G(x,y)dy$, note that $g(x+\alpha)=\int _{T^1}G(\psi(x,y))dy$. Hence, $(\ref{16'''})$ becomes in Fourier %%@
coefficients $k\neq 0$
\begin{multline}
\langle\mathbf{\ell}, (0,\hat{F}_1(k))(e^{2\pi ik\alpha}-1)\rangle+\hat{g}(k)(e^{2\pi ik2\alpha}-1)=\\
\langle \mathbf{\ell}_\gamma,C_0(0,\hat{F}_1(k))+[B_0+e^{2\pi ik\alpha}I]\hat{H}(k)\rangle
\end{multline}
A simple computation using $(\ref{1'})$ and $(\ref{2'})$ gives
\begin{equation}
\ell_1+\hat{g}(k)=\langle \mathbf{\ell}_\gamma,C_0(e^{2\pi ik\alpha}-1)^{-1}e_2+V+a(e^{2\pi ik\alpha}-1)^{-1}se_2\rangle
\end{equation}
Hence the equation $(\ref{16'''})$ has no  continuous solution $g$ by the Riemann Lebesgue Lemma. 
This implies necessarily that
$(\ref{3'})$ has no continuous solution $G$. This finishes this case.
 
Now let $p(x) =( x-1)^2 q_2( x)$. We may assume that ~\cite{N}
\begin{equation*}
L=\left( 
\begin{array}{cc}
A & 0 \\ 
C & B
\end{array}
\right)
\end{equation*}
where $1$ is not an eigenvalue of $B$, $A=\left( \begin{array}{cc}
1 & 0 \\ 
n & 1
\end{array}
\right)$ and $C\in M(2,\Z)$. The diffeomorphism $\varphi$ is given by 
\begin{equation}
\varphi(X,Y) =( AX +\mathbf{\delta}, CX +BY +F(x))
\label{B}
\end{equation}
where $X=(x,y)$, $Y=(z,w)$, $F(x) =(F_1(x) ,F_2(x))$ and $\mathbf{\delta}=(\alpha, \beta)$. Denoted by $a(X) =AX +\mathbf{\delta}$ the %%@
minimal affine transformation.

Suppose that $B$ is periodic with period $m$. Then
\begin{equation*}
\begin{array}{c}
\varphi ^m(X,Y) =(a^{m}(X), Y +C( m)X +\sum _{j=1}^{m}C( m-j-1) ( \alpha ,\beta ) +
\\ 
+\sum_{j=1}^{m}B^{m-j-1}F( x+( j-1) \alpha ))
\end{array}
\end{equation*}
Hence $( T^{4},\varphi ^{m}=\psi)$ is a simple free $T^2$-extension of the diffeomorphism $(T^2,a^m=\psi _0)$, given on covering $\R^2$ %%@
by $\psi_0(X)=a^{m}(X)$.
By ~\cite{P} $\psi$ is a minimal diffeomorphism if the equation $(\ref{7})$ has no continuous solution $f(X) =e^{2\pi i\left[ %%@
\left\langle \mathbf{\ell},X\right\rangle+G(X) \right]}$,  where $f_{\gamma }(X,Y) =e^{2\pi i\left\langle \mathbf{\ell }_{\gamma %%@
},Y\right\rangle }$ or equivalently the equation 

\begin{multline}
\left\langle \ell, a^{m}(X) -X\right\rangle +G(a^m(X)) -G(X)\\
=\left\langle \ell _{\gamma },C(m)X 
+\sum _{j=1}^{m}C( m-j-1)\delta+ \sum_{j=1}^{m}B^{m-j-1}F\left(x+( j-1)\alpha\right) \right\rangle
\label{18}
\end{multline}
has no continuous solution $G$. 
If
\begin{equation*}
\left\langle \ell, a^{m}(X) -X\right\rangle
\end{equation*} 
\begin{equation*}
\neq \left\langle \ell _{\gamma },C(m)X +\sum _{j=1}^{m}C( m-j-1)\delta \right\rangle
\end{equation*}
then one see that $(\ref{18})$ has no continuous solution $G$.
If 
\begin{equation*}
\left\langle \ell, a^{m}(X) -X\right\rangle
\end{equation*}
\begin{equation*}
= \left\langle \ell _{\gamma },(C(m)X +\sum _{j=1}^{m}C( m-j-1)\delta \right\rangle
\end{equation*}
the equation $(\ref{18})$ becomes
\begin{equation}
G( a^m(X)) -G(X)=\left\langle\ell _\gamma,\sum_{j=1}^{m}B^{m-j-1}F(x+(j-1)\alpha)\right\rangle
\end{equation}
or in Fourier coefficients 
\begin{equation*}
\hat{G}(k_1+nk_2,k_2)e^{2\pi i \langle a^m(0,0),(k_1,k_2)\rangle}-\hat{G}(k_1,k_2)
\end{equation*}
\begin{equation}
=\left\langle \ell_\gamma,\sum _{j=1}^m B^{m-j-1}\hat{F}(k_1)e^{2 \pi ik\alpha}\right\rangle
\label{a}
\end{equation}
If $k_2\neq 0$ then by the Riemann Lebesgue Lemma $\hat{G}(k_1,k_2)=0$
then the equation $(\ref{a})$ becomes

\begin{equation}
\hat{G}(k_1,0)=\left\langle \ell _\gamma,(e^{2 \pi ik\alpha}I-B)^{-1}\hat{F}(k_1)\right\rangle
\end{equation}
take $F$ as in $(\ref{2'})$ then the equation $\ref{18}$ has no continuous solution $G$.

If $B$ is not periodic, then 
\begin{equation*}
B=\left( 
\begin{array}{cc}
-1 & \mathbf{0} \\ 
s & -1
\end{array}
\right)
\end{equation*}
with $s\neq 0$. Consider the diffeomorphism as in $(\ref{B})$
\begin{multline}
\varphi(x,y)=(AX+\delta,CX+BY+F(x))\\
=(AX+\delta,C_1X-z+F_1(x),C_2X+sz-w+F_2(x))
\end{multline}
such that $\varphi _0(x,y,z)=(AX+\delta,C_1X-z+F_1(x))$ is minimal. Hence $(T^4,\psi=\varphi^2)$ is a simple free $T^1$-extension of the %%@
minimal diffeomorphism $(T^3,\psi _0)$, given on  covering $\R^3$ by $\psi=\varphi _0^{2}$. By ~\cite{P} $\psi$ is a minimal 
diffeomorphism if the equation $(\ref{7})$ has no continuous  solution $f(X)=e^{2\pi i[\langle \ell,X\rangle]+G(X)}$, where %%@
$f_\gamma(X,Y)=e^{2\pi i\ell_\gamma w}$ this equation is equivalent to the equation
\begin{multline}
\langle\ell,(a^2(X)-X,C_1(2)X+\delta(2)-F_1(x)+F_1(x+\alpha)\rangle+G(\psi_0(X))-G(X)\\
=\ell_\gamma[C_2(2)X-2sz+sF_1(x)-F_2(x)+F_2(x+\alpha)]
\end{multline}

The above equation  has no continuous solution $G$ because $s\neq 0$.

If $p(x)=(x-1)^3(x+1)$ then we assume that ~\cite{N} 
\begin{equation*}
L=\left( 
\begin{array}{cc}
A & 0 \\ 
C & -1
\end{array}
\right)
\end{equation*}
where $A\in GL(3,\Z)$ is unipotent, i.e., $A=I+N$. Consider the diffeomorphism $\varphi$ given by 
\begin{equation}
\varphi(X,w) =(AX+\mathbf{\alpha}, CX+-w+F(x))
\label{(-1)}
\end{equation}
where $X=(x,y,z)$, $\mathbf{\alpha}=(\alpha _1,\alpha _2,\alpha _3)$ and $F:\R\rightarrow \R$ is $\Z$-periodic smooth function. Then
\begin{multline}
\psi(X,w)=\varphi^2(X,w)=(A^2X+A\mathbf{\alpha}+\mathbf{\alpha},\\
[CA-C]X+w+C\mathbf{\alpha}-F(x)+F(x+\alpha _1))
\end{multline}

Thus $(T^4,\psi)$ is a simple free $T^1$-extension of the minimal affine transformation $(T^3,\psi _0)$, given on covering $\R^3$ by %%@
$\psi _0(X)=A^2X+A\mathbf{\alpha}+\mathbf{\alpha}$. Now by ~\cite{P} $\psi$ is minimal diffeomorphism if the equation $(\ref{7})$ has no %%@
continuous solution $f(X)=e^{2 \pi i[\langle \ell,X\rangle+G(X)]}$, where $f_\gamma(X,w)=e^{2 \pi i\ell _\gamma w}$ or equivalently the %%@
equation below has no continuous solution $G$.

\begin{multline}
\langle\ell, [A^2-I]X+A\mathbf{\alpha}+\mathbf{\alpha}\rangle+G(\psi _0(X))-G(X)=\\
\ell _\gamma[[CA-C]X+C\mathbf{\alpha}-F(x)+F(x+\alpha _1)]
\label{fin}
\end{multline}

If $\langle\ell, [A^2-I]X+A\mathbf{\alpha}\neq \ell _\gamma[[CA-C]X+C\mathbf{\alpha}]$ then one see that $(\ref{fin})$ has no continuous %%@
solution $G$. If $\langle\ell, [A^2-I]X+A\mathbf{\alpha}= \ell _\gamma[[CA-C]X+C\mathbf{\alpha}]$ the equation $(\ref{fin})$ becomes

\begin{equation}
G(\psi _0(X))-G(X)=\ell _\gamma [-F(x)+F(x+\alpha _1)]
\label{fin1}
\end{equation}
Now if $\alpha _1$ is a Liouville number and $F$ is a smooth function given by Fourier coefficients as in $(\ref{9})$ the equation %%@
$(\ref{fin1})$ has no continuous  solution $G$.
$\hfill \square$
\end{proof}
 
 We now give a sufficient condition for a smooth skew-product transformation of a torus to be smoothly  conjugate to an affine %%@
transformation. We present in the Appendix the definition of Diophantine vectors. We recall that we are restrict to smooth skew-product %%@
diffeomorphism $\varphi$ of the torus $T^n=T^p\times T^{n-p}$ given on the covering $\R^n$ by 
\begin{equation}
\varphi(X,Y)=(X+\alpha,CX+BY+F(X))
\label{fin2}
\end{equation} 
where $X\in\R^p$, $Y\in\R^{n-p}$ and $F:\R^p\longrightarrow\R^{n-p}$ is a smooth $\Z^p$-periodic function. We call $\alpha \in \R^p$ the %%@
translation vector of $\varphi$. To $\varphi$ there naturally corresponds an affine transformation $\varphi _0$ given on the covering %%@
$\R^n$ by 
$$
\varphi _0(X,Y)=(X+\alpha,CX+BY+\beta _1)
$$
where $\beta _1 =Proj(\beta)$, $Proj:\R^n\longrightarrow \ker(I-B)^t$ and $\beta=\int _{T^p}Fd\mu$, $\mu$ being the Haar measure of %%@
$T^p$. 
\begin{thm}
Every smooth skew-product diffeomorphism $\varphi$  of $T^n$ of type $(\ref{fin2})$ quasi-unipotent on homology whose translation vector %%@
is Diophantine is smoothly conjugate to its corresponding affine transformation. Moreover if $\varphi$ is minimal then $1$ is only %%@
eigenvalue of its linear part.
\end{thm}
\begin{proof} Since $\alpha$ translation vector is Diophantine and $B$ is quasi-\newline unipotent then cohomological equation that 
$$
F(x)=\beta _1+G(x+\alpha)-BG(x)
$$
has a smooth solution $G$. [~\cite{S}, Theorem 2.6]. Thus the diffeomorphism $h:T^n\longrightarrow T^n$ given on the covering $\R^n$ by 
$$
h(X,Y)=(x,Y+G(X))
$$
conjugate $\varphi$ with the affine transformation $\varphi _0$. If $\varphi$ is minimal then so is $\varphi _0$ and  $1$ is the only %%@
eigenvalue of $B$ ~\cite{H}.
$\hfill \square$
\end{proof}

\section{Appendix.}
In this section we recall the definition of Diophantine and Liouville numbers and we show that for every  nth root of unity there exist a %%@
"fast approximation" by iterations of a Liouville rotation of the circle.
\begin{definition}
Given $C>0$ and $r\geq 0$, we say $\alpha \in \R-\Q$ verifies a Diophantine condition of exponent $r$ and constant $C$ if and only if for %%@
all $q \in \Z$,  one has $||q\alpha||\geq C|q|^{-1-r}$.
\end{definition}

where  $||x||=\inf\{|x-p| |p\in\Z\}$. Notice that the inequality 
\begin{equation}
4s\leq|e^{2 \pi is}-1|\leq 2\pi s
\label{ap1}
\end{equation}
with $s\in [0,1]$ implies that the orbit of $1$ by the rotation $R_{\alpha}$ is a bad approximation of number $1$ in the sense that
\begin{equation}
|e^{2 \pi i q\alpha}-1|\geq 4C|q|^{-1-r}
\end{equation}
for all $q\in \Z$.

The irrational algebraic numbers are examples of Diophantine numbers.

More generally we say that a vector $\alpha \in \R^n$, $n>1$ is a Diophantine vector if there are constants $r>0$ and $C>0$ such that 

\begin{equation}
||\langle k,\alpha \rangle||\geq C|k|^{-r}
\label{ap2}
\end{equation}
for every $k\in \Z^n$.
 
We say that an irrational number is Liouville if it is not Diophantine. In this case the orbit of $1$ by the rotation $R_\alpha$ has a %%@
good approximation of the number $1$ in the sense that: there exists a sequence $\{q_j\}_{j\in \N}\subset\Z$ such that %%@
$|\alpha-\frac{p_j}{q_j}|<C|q_j|^{-(j+1)}$ por $(\ref{ap1})$.

The number $\alpha=\sum^{\infty} _{k=1}q^{-k!}$ is a Liouville number.

The following proposition shows that for every root of unity $\xi _n$ there exists a rotation $R_\alpha$ whose orbit by $1$ has fast %%@
approximation to $\xi _n$.

\begin{prop} Given a family $\{e^{2\pi i \frac{p_s}{q_s}}\}^m_{s=1}$ of roots of unity there exits a Liouville number $\alpha$ and a %%@
family $\{\{k^s_j\}_{j\in \N}\}^m_{s=1}$ of sequences such that  
\begin{equation}
||k^s_j\alpha-\frac{p_s}{q_s}||<\frac{2}{(k^s_j)^j}
\label{a1} 
\end{equation}
para todo $j\in \N$ and $s\in \{1,2..,m\}$

\end{prop}
\begin{proof} We may assume that $q_s>0$ y  $0<p_s<q_s$. We consider the following Liouville number 
$$
\alpha=\sum^\infty _{k=1}\frac{1}{q^{k!}}
$$
where $q=q_1q_2....q_m$. Now consider the sequence $k^s_j=p_sq^{j!}_1...(q_s)^{j!-1}...q^{j!}_m$, con $j\in \N$ y $s\in \{1,2,...,m\}$. %%@
We will show that this sequence satisfies $(\ref{a1})$. In fact, it is easy to see that

$||k^s_j\alpha-\frac{p_s}{q_s}||=\inf_{\ell \in \Z}|k^s_j\alpha-\frac{p_s}{q_s}+\ell|\leq\sum^\infty _{k=j+1}\frac{k^s_j}{q^{k!}}$.

Thus
\begin{multline}
||k^s_j\alpha-\frac{p_s}{q_s}||\leq \frac{1}{(k^s_j)^j}\sum^\infty %%@
_{k=j+1}\frac{(k^s_j)^{j+1}}{q^{k!}}\\\leq\frac{1}{(k^s_j)^j}\sum^\infty _{k=j+1}\frac{q^{(j+1)!}}{q^{k!}}<\frac{1}{(k^s_j)^j}\sum^\infty %%@
_{k=0}\frac{1}{q^k}\leq\frac{2}{(k^s_j)^j}
\end{multline}
$\hfill \square$
\end{proof}
\section*{}

\end{document}